\documentclass[a4paper, 11pt, twoside]{article}
\usepackage{mfpaperstuff}
\usepackage{tikz-3dplot}
\usetikzlibrary{shapes}
\begin{document}
\medmuskip=3mu plus 1mu minus 1mu
\newcommand\ci[2]{#1{\circ}#2}
\newcommand\ltup[2]{\left(\left.#1\ \right|\ \smash{#2}\right)}
\newcommand\rtup[2]{\left(\smash{#1}\ \left|\ #2\right.\right)}
\newcommand\lmst[2]{\left[\left.#1\ \right|\ \smash{#2}\right]}
\newcommand\rmst[2]{\left[\smash{#1}\ \left|\ #2\right.\right]}
\newcommand\plu[2]{#1\boxplus#2}
\newcommand\weyl{W_}
\newcommand{\thmcite}[2]{\textup{\textbf{\cite[#2]{#1}}}\ }
\newcommand\lra\longrightarrow
\newcommand\cornoarg[1]{\mathrm{cor}_{#1}}
\newcommand\weinoarg[1]{\mathrm{wt}_{#1}}
\newcommand\cor[2]{\cornoarg{#1}#2}
\newcommand\wei[2]{\weinoarg{#1}#2}
\newcommand\quo[2]{#2^{(#1)}}
\newcommand\qou[1]{\mathrm{quo}_{#1}}
\newcommand\be[1]{\mathcal{B}^{#1}}
\newcommand\ber[2]{\calb_{#1}^{#2}}
\newcommand\gc[2]{$[#1{:}#2]$-core}
\newcommand\gcc[3]{[$#1{:}#2{:}#3$]-core}
\newcommand\cg[2]{\mathcal{C}_{#1:#2}}
\newcommand\cores[1]{\calc_{#1}}
\newcommand\cnd{(\textasteriskcentered)}
\newcommand\cne{(\textdagger)}
\newcommand\zsz{\bbz/s\bbz}
\newcommand\ztz{\bbz/t\bbz}
\newcommand\es[2]{\Gamma_{#1}#2}
\newcommand\twq[1]{\mathrm Q_s^t#1}

\Yvcentermath1

\title{A generalisation of core partitions}
\toptitle

\begin{abstract}
Suppose $s$ and $t$ are coprime natural numbers. A theorem of Olsson says that the $t$-core of an $s$-core partition is again an $s$-core. We generalise this theorem, showing that the $s$-weight of the $t$-core of a partition $\la$ is at most the $s$-weight of $\la$. Then we consider the set $\cg st$ of partitions for which equality holds, which we call \gc sts; this set has interesting structure, and we expect that it will be the subject of future study. We show that the set of \gc sts is a union of finitely many orbits for an action of a Coxeter group of type $\tilde A_{s-1}\times\tilde A_{t-1}$ on the set of partitions. We also consider the problem of constructing an \gc st with specified $s$-core and $t$-core.
\end{abstract}

\section{Introduction}

In this paper we study the combinatorics of integer partitions. If $s$ is a natural number, an \emph{$s$-core} (often referred to in the literature as an \emph{$s$-core partition}) is a partition with no rim hooks of length $s$. The $s$-core of an arbitrary partition $\la$ is obtained by repeatedly removing rim $s$-hooks from $\la$, and the $s$-weight of $\la$ is the number of rim hooks removed. The set of all $s$-cores, which we denote $\cores s$, has geometric structure related to type A alcove geometry, and has applications in representation theory and number theory.

Now suppose $t$ is another natural number which is prime to $s$. A recent trend in partition combinatorics has been to compare $s$-cores and $t$-cores. Anderson \cite{and} enumerated \emph{$(s,t)$-cores}, i.e.\ partitions which are both $s$- and $t$-cores, and Fishel and Vazirani \cite{fv} explored the connection between $(s,t)$-cores and the associated alcove geometry. Several authors \cite{k,os,vand,mfcores} have studied the properties of the largest $(s,t)$-core, which is denoted $\kappa_{s,t}$. The present author explored another avenue in \cite{mfcores}, considering the $t$-core of an arbitrary $s$-core; by a theorem of Olsson, the $t$-core of an $s$-core is again an $s$-core, so we have a natural map from $s$-cores to $(s,t)$-cores. Exploiting natural actions of the Coxeter group $\weyl s$ of type $\tilde A_{s-1}$ on $\cores s$ reveals interesting symmetry in this map.

In the present paper we generalise Olsson's theorem, showing that replacing any partition with its $t$-core does not increase its $s$-weight. We define a an \emph{\gc st} to be a partition for which equality holds in this statement, and explore the family $\cg st$ of \gc sts, which plays a kind of dual role to the family of $(s,t)$-cores. We show that $\cg st$ is a union of orbits for a certain action of $\weyl s\times\weyl t$, with each orbit containing a unique $(s,t)$-core. We then consider the problem of constructing a partition with a given $s$-core $\sigma$ and a given $t$-core $\tau$; we show that if the $t$-core of $\sigma$ coincides with the $s$-core of $\tau$, then there is a unique \gc st with $s$-core $\sigma$ and $t$-core $\tau$, and we give a simple method for constructing this partition. This leads to an alternative characterisation of an \gc st as a partition which is uniquely determined by its size, its $s$-core and its $t$-core. Finally we consider the orbit of $\weyl s\times\weyl t$ containing $\kappa_{s,t}$, showing that this is naturally in bijection with $\cores s\times\cores t$.

In the next section we recall basic definitions and simple results, largely to fix conventions and notation. In \cref{groupsec} we define the group $\weyl s$ and study its actions on integers and partitions; some of this material has not appeared in this form before. In \cref{gcsec,sumsec,kappasec} we prove our main results. We finish with some brief comments in \cref{rmksec}.

\section{Basic definitions}\label{defsec}

\subsection{Conventions and notation}

In this section we set out some basic conventions that we use throughout the paper. As usual, $\bbn$ denotes the set of positive integers, and $\bbn_0$ the set of non-negative integers. We shall often consider the set $\zsz$, where $s\in\bbn$, and we use the formal convention that
\[
\zsz=\lset{a+s\bbz}{a\in\bbz},
\]
where as usual
\[
a+s\bbz=\lset{a+sb}{b\in\bbz}
\]
for an integer $a$. We adopt the standard convention that $(a+s\bbz)+b=(a+b)+s\bbz$ for any $a,b\in\bbz$, but we adopt an unusual convention for multiplication, namely that
\[
(a+s\bbz)b=(ab)+s\bbz.
\]

If $X$ is a set, a \emph{$\zsz$-tuple} of elements of $X$ is simply a function $i\mapsto x_i$ from $\zsz$ to $X$, and we may write such a tuple in the form $\ltup{x_i}{i\in\zsz}$. A \emph{multiset} of elements of $X$ is a subset of $X$ with possibly repeated elements (i.e.\ a function from $X$ to $\bbn_0$). We write a multiset by writing the elements (with multiplicity) surrounded by square brackets. Given a $\zsz$-tuple $\ltup{x_i}{i\in\zsz}$, we denote the associated multiset $\lmst{x_i}{i\in\zsz}$.

\subsection{Partitions}

In this paper, a \emph{partition} is an infinite weakly decreasing sequence $\la=(\la_1,\la_2,\dots)$ of non-negative integers with finite sum; we write $|\la|$ for this sum. When writing partitions, we omit zeroes and group together equal parts with a superscript, so the partition $(4,3,3,1,1,0,0,\dots)$ is written as $(4,3^2,1^2)$. The partition $(0,0,\dots)$ is denoted $\varnothing$, and the set of all partitions is denoted~$\calp$.

The \emph{Young diagram} of a partition $\la$ is the set of all pairs $(i,j)\in\bbn^2$ for which $j\ls\la_i$; we often identify $\la$ with its Young diagram, so for example we may write $\la\subseteq\mu$ to mean that $\la_i\ls\mu_i$ for all $i$. The \emph{rim} of $\la$ is the set of all $(i,j)\in\la$ such that $(i+1,j+1)\notin\la$.
We draw the Young diagram as an array of boxes in the plane; for example, the Young diagram of $(4,3^2,1^2)$ (with the elements of the rim marked $\times$) is as follows:
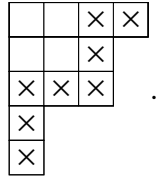
\newcommand\rti{\raisebox{1pt}{$\times$}}
\[
\young(\ \ \rti\rti,\ \ \rti,\rti\rti\rti,\rti,\rti).
\]

\subsection{Cores}

Suppose $\la$ is a partition and $s\in\bbn$. A \emph{rim $s$-hook} of $\la$ is a subset of the rim of $\la$ of size exactly $s$ which is connected (i.e.\ comprises consecutive elements of the rim) and which is removable in the sense that it can be removed to leave the Young diagram of a partition. For example, in the example $\la=(4,3^2,1^2)$ above, the set $\{(2,3),(3,3),(3,2)\}$ is a rim $3$-hook.

The \emph{$s$-core} of $\la$ is the partition obtained from $\la$ by repeatedly removing rim $s$-hooks until none remain. for example, the $3$-core of the partition $\la=(4,3^2,1^2)$ above is $(4,2)$. It is well known (and follows from \cref{basicbeta} below) that the $s$-core of $\la$ is independent of the choice of rim hook removed at each stage. The \emph{$s$-weight} of $\la$ is the number of rim $s$-hooks removed to reach the $s$-core of $\la$, i.e.\ $\frac1s(|\la|-|\cor s\la|)$. We write $\cor s\la$ for the $s$-core of $\la$, and $\wei s\la$ for the $s$-weight. Note that these definitions remain valid in the case $s=1$, although this case is seldom considered in the literature. In this case, we have $\cor1\la=\varnothing$ for any $\la$, and hence $\wei1\la=|\la|$.

We say that $\la$ is an $s$-core if $\cor s\la=\la$ (or equivalently if $\wei s\la=0$), and we write $\cores s$ for the set of all $s$-cores. This set has been studied at length; it enjoys a rich geometric structure, and has applications in representation theory and number theory.

A trend in recent years has been to compare $s$-cores and $t$-cores, where $s$ and $t$ are distinct positive integers. We define an \emph{$(s,t)$-core} to be a partition which is both an $s$-core and a $t$-core. It is known that there are finitely many $(s,t)$-cores if and only if $s$ and $t$ are coprime; in this case, the number of $(s,t)$-cores is precisely $\mfrac1{s+t}\mbinom{s+t}s$ \cite[Theorems 1 \& 3]{and}.

\subsection{Beta-sets}

A useful way to understand partitions, and in particular $s$-cores, is via beta-sets. Given a partition $\la$ and an integer $r$, we define the \emph{beta-set} $\ber r\la$ to be the infinite set of integers
\[
\ber r\la=\lset{\la_i-i+r}{i\in\bbn}.
\]
We shall mostly consider the beta-set $\ber0\la$, which we denote simply $\be\la$.

Note that $\ber r\la$ is bounded above and its complement in $\bbz$ is bounded below. Conversely, if we are given a subset $\calb$ of $\bbz$ which is bounded above and whose complement is bounded below, then we have $\calb=\ber r\la$ for some (uniquely defined) partition $\la$ and integer $r$: we let $r$ be the number of positive integers in $\calb$ minus the number of non-positive integers not in $\calb$; then, writing the elements of $\calb$ as $b_1>b_2>\dots$ and setting $\la_i=b_i+i-r$ for each $i$, we have a partition $\la$, and $\calb=\ber r\la$.

Later we will need the following simple lemma.

\begin{lemma}\label{becont}
Suppose $\la$ and $\mu$ are partitions and $r\in\bbz$. If $\ber r\la\subseteq\ber r\mu$, then $\la=\mu$.
\end{lemma}

\begin{pf}
Choose $N$ sufficiently large that $\la_N=\mu_N=0$. Then $\ber r\la$ and $\ber r\mu$ both contain all integers less than or equal to $r-N$. Moreover, $\be\la$ contains exactly $N-1$ elements greater than $r-N$, namely $\la_1-1+r,\dots,\la_{N-1}-(N-1)+r$, and similarly $\be\mu$ contains exactly $N-1$ elements greater than $r-N$. Since $\ber r\la\subseteq\ber r\mu$, we get $\ber r\la=\ber r\mu$, and hence $\la=\mu$.
\end{pf}

The main advantage of beta-sets is the easy identification and classification of $s$-cores. The following lemma is due to James, and is the main motivation for introducing beta-sets.

\begin{lemma}\label{basicbeta}
Suppose $\la$ is a partition, $r\in\bbz$ and $s\in\bbn$. Then $\la$ is an $s$-core if and only if there is no $b\in\ber r\la$ such that $b-s\notin\ber r\la$; moreover, the beta-set $\ber r{\cor s\la}$ can be obtained from $\ber r\la$ by repeatedly replacing integers $b$ with $b-s$ as far as possible. 
\end{lemma}

\begin{eg}
Take $\la=(4,3^2,1^2)$. Then $\be\la=\{3,1,0,-3,-4,-6,-7,\dots\}$. Taking $s=5$, we see that there are two integers $a\in\be\la$ such that $a-s\notin\be\la$, namely $3$ and $0$. Replacing $3$ with $-2$ or $0$ with $-5$ corresponds to removing a rim $5$-hook from $\la$, and making both of these replacements yields the set $\{1,-2,-3,-4,\dots\}$, which is the beta-set of the partition $(2)=\cor5\la$.
\end{eg}

A very helpful way to visualise a beta-set of a partition is via James's \emph{abacus}. Take an abacus with $s$ infinite vertical runners, numbered $0,\dots,s-1$ from left to right, and mark positions on the runners labelled by the integers, such that position $x$ is immediately below position $x-s$ for all $x$, and position $x$ is immediately to the right of position $x-1$ whenever $s\nmid x$. For example, when $s=3$ the positions are marked as follows.
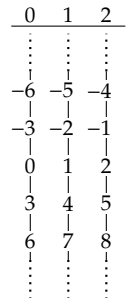
\[
\begin{tikzpicture}[scale=.5,every node/.style={fill=white,inner sep=1pt}]
\foreach\x in{0,1,2}{
\draw(\x,-2.5)--++(0,5);
\draw[dotted,thick](\x,2.5)--++(0,1);
\draw[dotted,thick](\x,-2.5)--++(0,-1);
}
\draw(-.5,3.7)--++(3,0);
\foreach\x in{0,1,2}\node at(\x,4){\scriptsize$\x$};
\node at(0,2){\scriptsize$\llap{$-$}6$};
\node at(1,2){\scriptsize$\llap{$-$}5$};
\node at(2,2){\scriptsize$\llap{$-$}4$};
\node at(0,1){\scriptsize$\llap{$-$}3$};
\node at(1,1){\scriptsize$\llap{$-$}2$};
\node at(2,1){\scriptsize$\llap{$-$}1$};
\node at(0,0){\scriptsize$0$};
\node at(1,0){\scriptsize$1$};
\node at(2,0){\scriptsize$2$};
\node at(0,-1){\scriptsize$3$};
\node at(1,-1){\scriptsize$4$};
\node at(2,-1){\scriptsize$5$};
\node at(0,-2){\scriptsize$6$};
\node at(1,-2){\scriptsize$7$};
\node at(2,-2){\scriptsize$8$};
\end{tikzpicture}
\]
Now given a partition $\la$ and an integer $r$, place a bead on the abacus at position $x$ for each $x\in\ber r\la$. The resulting configuration is called an \emph{abacus display} for $\la$. \Cref{basicbeta} can now be stated as follows: $\la$ is an $s$-core if and only if every bead in an abacus display for $\la$ has a bead immediately above it; if $\la$ is not an $s$-core, an abacus display for the $s$-core of $\la$ can be obtained by sliding the beads up their runners as far as possible.

\begin{eg}
Taking $\la=(4,3^2,1^2)$ and $s=5$ as in the previous example, we obtain the following abacus displays for $\la$ and $\cor5\la$. (Note that we typically omit the labels of runners and positions when drawing abacus displays.)
\[
\abacus(vvvvv,bbbbb,bbbbb,nbbnn,bbnbn,nnnnn,nnnnn;vvvvv)
\qquad
\abacus(vvvvv,bbbbb,bbbbb,bbbbn,nbnnn,nnnnn,nnnnn;vvvvv)
\]
\end{eg}

\cref{basicbeta} implies the following two results.

\begin{cory}\label{weightbeta}
Suppose $\la$ is a partition, $r\in\bbz$ and $s\in\bbn$. Then $\wei s\la$ equals the number of pairs $(x,l)\in\bbz\times\bbn$ such that $x\in\ber r\la\notni x-ls$.
\end{cory}

\begin{cory}\label{scoremscore}
Suppose $r,s\in\bbn$ and $\la\in\calp$. Then $\cor s\la=\cor s(\cor{rs}\la)$. In particular, if $\la$ is an $s$-core then $\la$ is an $rs$-core.
\end{cory}

\subsection{$s$-sets}

Now we define the $s$-set of a partition $\la$; this is a set of $s$ integers which provides a useful encoding of the $s$-core of $\la$.

Suppose $\la$ is an $s$-core. Given $i\in\zsz$, we define $\es i\la$ to be the smallest integer in $i$ which is not in $\be\la$. Then by \cref{basicbeta} we have
\[
\be\la\cap i=\lset{\es i\la-ks}{k\in\bbn}.
\]
Following \cite{mfcores}, we refer to the set $\lset{\es i\la}{i\in\zsz}$ as the \emph{$s$-set} of $\la$; this set consists of $s$ integers which are pairwise incongruent modulo $s$, and which sum to $\frac12s(s-1)$. Conversely, any set of $s$ integers with these properties is the $s$-set of some $s$-core. In terms of the abacus, the $s$-set consists of those unoccupied positions $x$ such that $x-s$ is occupied, in the abacus display for $\la$ with $r=0$.

Now suppose $\la$ is an arbitrary partition. We define $\es i\la=\es i{(\cor s\la)}$ for each $i$, and refer to the $s$-set of $\la$ meaning the $s$-set of $\cor s\la$.

\subsection{Quotients and the abacus}\label{quotsec}

Next we define the $s$-quotient of a partition $\la$. Given $j\in\zsz$, write $j=i+s\bbz$ for some integer $i$, and consider the set $\be\la\cap j$. By subtracting $i$ from each of the integers in this set and dividing by $s$, we obtain a set of integers which is bounded above and whose complement in $\bbz$ is bounded below. This set is therefore a beta-set of some partition, which we denote $\quo j\la$. It is very easy to see that this partition is independent of the choice of $i$.  We define the \emph{$s$-quotient} $\qou s\la$ to be the $\zsz$-tuple $\rtup{\quo j\la}{\smash j\in\zsz}$.

\begin{eg}
Suppose $s=3$ and $\la=(5,4^2,3,2,1^2)$. Then
\[
\be\la=\{4,2,1,-1,-3,-5,-6,-8,-9,\dots\},
\]
so
\[
\be\la\cap(1+3\bbz)=\{4,1,-5,-8,-11,\dots\}.
\]
Subtracting $1$ from each element and dividing by $3$, we obtain the set $\{1,0,-2,-3,-4,\dots\}$. This is the beta-set $\ber1{(1^2)}$, so we have $\quo{1+3\bbz}\la=(1^2)$. In a similar way we find that $\quo{3\bbz}\la=\varnothing$ and $\quo{2+3\bbz}\la=(2^2)$.
\end{eg}

The $s$-quotient of a partition can easily be visualised on the abacus. Taking the abacus display for $\la$ with $s$ runners and with $r=0$, examine runner $i$ in isolation; this can be regarded as a $1$-runner abacus display for a partition, and this partition is $\quo{i+s\bbz}\la$. In other words, $\quo{i+s\bbz}\la_l$ equals the number of unoccupied positions above the $l$th lowest bead on runner $i$.

\begin{eg}
Taking $\la=(5,4^2,3,2,1^2)$ and $s=3$ as in the previous example, we obtain the following abacus display, from which we can see that the $3$-quotient of $\la$ is as given above.
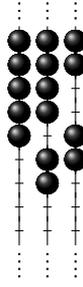
\[
\abacus(vvv,bbb,bbb,bbn,bbn,bnb,nbb,nbn,nnn,nnn;vvv)
\]
\end{eg}

Note that a partition $\la$ is determined by its $s$-core and $s$-quotient. To see this, let $\sigma=\cor s\la$, and for each $j\in\zsz$ consider the set $\be\sigma\cap j$. By \cref{basicbeta} we have
\begin{align*}
\be\sigma\cap j&=\rset{\es j\sigma-is}{i\in\bbn}.\\
\intertext{So if we let $\upsilon=\quo j\mu$, then}
\be\la\cap j&=\rset{\es j\sigma+\upsilon_is-is}{i\in\bbn}.
\end{align*}
Applying this for all $j$, we see that $\be\la$ (and hence $\la$) is determined by $\sigma$ and $\qou s\la$. Moreover, we can see that for any $s$-core $\sigma$ and any $\zsz$-tuple $\upsilon=\rtup{\upsilon^{(j)}}{\smash j\in\zsz}$ of partitions, there is a partition with $s$-core $\sigma$ and $s$-quotient $\upsilon$.

Quotients of partitions will prove useful below. Two important properties are given in the following \lcnamecref{quotprops}.

\begin{lemma}\label{quotprops}
Suppose $\la\in\calp$ and $r,s\in\bbn$.
\begin{enumerate}
\item\label{wc1}
$\wei{rs}\la=\sum_{j\in\zsz}\wei r{\quo j\la}$; in particular, $\la$ is an $rs$-core if and only if each component of $\qou s\la$ is an $r$-core.
\item\label{wc2}
For each $j\in\zsz$, we have $\quo j{(\cor{rs}\la)}=\cor{r}{(\quo j\la)}$.
\end{enumerate}
\end{lemma}

\begin{pf}
Both statements follow from \cref{basicbeta}: removing a rim $rs$-hook from $\la$ corresponds to reducing some element of $\be\la$ by $rs$; this in turn corresponds to reducing an element of the beta-set of one component of $\qou s\la$ by $r$.
\end{pf}

\section{The affine symmetric group}\label{groupsec}

In this section we assume $s\gs2$. Let $\weyl s$ denote the Coxeter group of type $\tilde A_{s-1}$; this has generators $w_i$ for $i\in\zsz$, and relations
\begin{alignat*}2
w_i^2&=1\qquad&&\text{for each }i,\\
w_iw_j&=w_jw_i&&\text{if }j\neq i\pm1,\\
w_iw_jw_i&=w_jw_iw_j&\qquad&\text{if }j=i+1\neq i-1.
\end{alignat*}
Now suppose $t$ is another positive integer which is prime to $s$, and set $\ci st=\frac12(s-1)(t-1)$. We define the \emph{level $t$ action} of $\weyl s$ on $\bbz$ by
\[
w_in=\begin{cases}
n+t&(n\in(i-1)t-\ci st)\\
n-t&(n\in it-\ci st)\\
n&(\text{otherwise})
\end{cases}\tag*{for each $i\in\zsz$.}
\]
(Recall from \cref{defsec} our convention: if $i=a+s\bbz$, then $it$ means $at+s\bbz$.)

This action is faithful for every $t$; in fact, the image of the level $1$ action is often taken as a concrete definition of $\weyl s$.

It is easy to see that if $\calb$ is a subset of $\bbz$ which is bounded above and whose complement in $\bbz$ is bounded below, then the same is true of $w\calb$ for any $w\in\weyl s$; moreover, the number of non-negative elements minus the number of negative non-elements is the same in $\calb$ and $w\calb$. Hence we have an action of $\weyl s$ on $\calp$, given by
\[
\be{w_i\la}=w_i\be\la.
\]
We refer to this action as the \emph{level $t$ action of $\weyl s$ on $\calp$}, and we refer to an orbit under this action as a \emph{level $t$ orbit}.

\begin{rmk}
The level $1$ action of $\weyl s$ on $\calp$ is well known, and was first addressed by Lascoux \cite{l}. In \cite{mfcores}, the author introduced the level $t$ action of $\weyl s$ on $\bbz$ and on the set of $s$-cores (which is a union of orbits for the action on $\calp$), but with a slight difference from the definition above, in that the terms $-\ci st$ do not appear in the definition in \cite{mfcores}. This makes little practical difference, since the two versions of the action are equivalent to each other via a diagram automorphism of $\weyl s$. However, we prefer the version above in this paper; although slightly more complicated to define, it turns out to be more helpful, as we shall see in \cref{bijec}. It also respects conjugation of partitions, in the sense that $(w_i\la)'=w_{-i}(\la')$ for any $i$ and any $\la$, where $\la'$ denotes the conjugate (i.e.\ transpose) partition.
\end{rmk}

\begin{eg}
Take $s=2$, $t=3$, so that $\ci st=1$. Let $\la=(4,1)$, which has beta-set $\be\la=\{3,-1,-3,-4,\dots\}$. To calculate $w_{1+2\bbz}\la$, we add $3$ to every odd element of this beta-set and subtract $3$ from every even element; we obtain $\{6,2,0,-2,-4,-6,-7,-8,\dots\}=\be{(7,4,3,2,1)}$, so $w_{1+2\bbz}\la=(7,4,3,2,1)$. in a similar way, we calculate $w_{2\bbz}\la=(1^2)$. Note that $w_{1+2\bbz}\la$ is obtained by adding four rim $3$-hooks to $\la$, while $w_{2\bbz}\la$ is obtained by removing a rim $3$-hook.  In general, the effect of $w_i$ acting on a partition $\la$ at level $t$ can be described in terms of simultaneously adding and removing rim $t$-hooks; we leave the reader to work out the details.
\end{eg}

We now give some invariants of the level $t$ action of $\weyl s$. In \cref{samewquot} we shall use these to give an explicit criterion for when two partitions lie in the same level $t$ orbit.

\begin{lemma}\label{orbitprops}
Suppose $\la\in\calp$ and $w\in\weyl s$, and define $w\la$ using the level $t$ action.  Then:
\begin{enumerate}
\item\label{op1}
$\cor t{(w\la)}=\cor t\la$;
\item\label{op1a}
$\qou s{(w\la)}$ is the same as $\qou s\la$ with the components re-ordered;
\item\label{op2}
$\wei s{(w\la)}=\wei s\la$;
\item\label{corcomm}
$\cor s{(w\la)}=w{(\cor s\la)}$.
\end{enumerate}
\end{lemma}

\begin{pf}
Since the relations occurring in all four parts are transitive, we may assume $w=w_i$ for $i\in\zsz$. Write $j=it-\ci st$.
\begin{enumerate}
\item
We obtain $\be{w_i\la}$ from $\be\la$ by adding $t$ to every element of $\be\la\cap(j-t)$ and subtracting $t$ from every element of $\be\la\cap j$.  But we may as well ignore those pairs of integers $b,b-t$ which both lie in $\be\la$ and for which $b\in j$.  Since all but finitely many elements $b\in\be\la$ satisfy $b-t,b+t\in\be\la$, this means that we can get from $\be\la$ to $\be{w_i\la}$ in a finite sequence of moves, where each move is either increasing an element by $t$ or decreasing an element by $t$.  In other words, we can get from $\la$ to $w_i\la$ by adding and removing finitely many $t$-hooks.  So $\la$ and $w_i\la$ have the same $t$-core.
\item
For any $l\in\zsz\setminus\{j,j-t\}$, we have $\be{w_i\la}\cap l=\be{\la}\cap l$, so that $\quo l{(w_i\la)}=\quo l\la$. On the other hand, we have
\begin{align*}
\be{w_i\la}\cap(j-t)&=\lset{b-t}{b\in\be\la\cap j}\\
\be{w_i\la}\cap j&=\lset{b+t}{b\in\be\la\cap(j-t)},
\end{align*}
which gives $\quo{j-t}{(w_i\la)}=\quo j\la$ and $\quo j{(w_i\la)}=\quo{j-t}\la$.
\item
This follows from (\ref{op1a}) and \cref{quotprops}(\ref{wc1}) (taking $r=1$ in that \lcnamecref{quotprops}).
\item
From the definition of the $s$-quotient, the largest element of $\be\la\cap l$ is $s\quo l\la_1+\es l\la-s$. Taking $l=j$ and applying $w_i$, the largest element of $\be{(w_i\la)}\cap(j-t)$ is $s\quo j\la_1+\es j\la-s-t$, and this must equal $s\quo{j-t}{(w_i\la)}_1+\es{j-t}{(w_i\la)}-s$. From the proof of part (\ref{op1a}) we have $\quo{j-t}{(w_i\la)}=\quo j\la$, and we deduce that $\es{j-t}{(w_i\la)}=\es j\la-t$. Similarly $\es j{(w_i\la)}=\es{j-t}\la+t$, and obviously $\es l{(w_i\la)}=\es l\la$ for $l\neq j,j-t$.

In particular, the $s$-set of $w_i\la$ is determined by the $s$-set of $\la$. Since $\la$ and $\cor s\la$ have the same $s$-set, so do $w_i\la$ and $w_i(\cor s\la)$. By (\ref{op2}) $w_i(\cor s\la)$ is an $s$-core, and since there is a unique $s$-core with a given $s$-set, we therefore have $w_i(\cor s\la)=\cor s(w_i\la)$.
\qedhere
\end{enumerate}
\end{pf}

Next we give a criterion for determining when two partitions lie in the same level $t$ orbit; to do this, we shall need to cite a result from \cite{mfcores} which gives a condition for two $s$-cores to lie in the same level $t$ orbit. (Note that although a slightly different level $t$ action is used in that paper, it differs from our action only by an automorphism of $\weyl s$, and so the orbits for the two actions are the same.)

\begin{propnc}{mfcores}{Proposition 4.1 \& Corollary 4.5}\label{4145}
Suppose $\la$ and $\mu$ are $s$-cores, and that the multisets
\[
\lmst{\es i\la+t\bbz}{i\in\zsz},\qquad\lmst{\es i\smash\mu+t\bbz}{i\in\zsz}
\]
are equal. Then $\cor t\la=\cor t\mu$, and $\la$ and $\mu$ lie in the same level $t$ orbit.
\end{propnc}

For our more general result, we make a definition which combines the $s$-quotient of a partition with its $s$-set modulo $t$.

Suppose $\la\in\calp$, with $s$-set $\lset{\es i\la}{i\in\zsz}$ and $s$-quotient $\rtup{\quo i\la}{i\in\zsz}$. We define the \emph{$t$-weighted $s$-quotient} of $\la$ to be the multiset
\[
\twq\la=\rmst{(\es i\la+t\bbz,\quo i\la)}{i\in\zsz}.
\]
of elements of $\ztz\times\calp$.

\begin{eg}
Take $s=4$ and $\la=(10,8,7,5,2,1^4)$. Then we have
\begin{alignat*}4
\es{4\bbz}\la&=0,\qquad&\es{1+4\bbz}\la&=9,&\qquad\es{2+4\bbz}\la&=2,&\qquad\es{3+4\bbz}\la&=-5,\\
\quo{4\bbz}\la&=(2),&\quo{1+4\bbz}\la&=(1),&\quo{2+4\bbz}\la&=(2),&\quo{3+4\bbz}\la&=(1).
\end{alignat*}
So the $7$-weighted $4$-quotient of $\la$ is
\[
[(7\bbz,(2)),(2+7\bbz,(1)),(2+7\bbz,(2)),(2+7\bbz,(1))].
\]
\end{eg}

\begin{propn}\label{samewquot}
Suppose $\la,\mu\in\calp$, and $s$ and $t$ are coprime positive integers. Then $\la$ and $\mu$ lie in the same level $t$ orbit of $\weyl s$ if and only if they have the same $t$-weighted $s$-quotient.
\end{propn}

\begin{pf}
For the `only if' part, we may assume that $\mu=w_i\la$ for $i\in\zsz$, and we write $j=it-\ci st$. From the proof of \cref{orbitprops} we have
\[
\quo l\mu=\begin{cases}
\quo j\la&(l=j-t)\\
\quo{j-t}\la&(l=j)\\
\quo l\la&(\text{otherwise}),
\end{cases}
\qquad
\es l\mu=\begin{cases}
\es j\la-t&(l=j-t)\\
\es{j-t}\la+t&(l=j)\\
\es j\la&(\text{otherwise}).
\end{cases}
\]
So $\la$ and $\mu$ have the same $t$-weighted $s$-quotient.

For the `if' part, assume that $\la$ and $\mu$ have the same $t$-weighted $s$-quotient. Since $\cor s\la$ and $\la$ have the same $s$-set (by definition) and the $s$-quotient of $\cor s\la$ has all components equals to $\varnothing$, we see that $\cor s\la$ and $\cor s\mu$ also have the same $t$-weighted $s$-quotient. So by \cref{4145} $\cor s\la$ and $\cor s\mu$ have the same $t$-core $\xi$, and lie in the same level $t$ orbit as $\xi$; that is, there are $w,x\in\weyl s$ such that $w\cor s\la=\xi=x\cor s\mu$. By \cref{orbitprops}(\ref{corcomm}) we have $\cor s{(w\la)}=\xi=\cor s{(x\mu)}$, so by replacing $\la$ and $\mu$ with $w\la$ and $x\mu$ (and using the `only if' part above), we may assume that $\la$ and $\mu$ both have $s$-core $\xi$, with $\xi$ being an $(s,t)$-core. So we have $\es i\la=\es i\mu=\es i\xi$ for every $i$. Now given any $r\in\ztz$, let
\[
\calx_r=\lset{i\in\zsz}{\es i\xi\in r}.
\]
Then the fact that $\la$ and $\mu$ have the same $t$-weighted $s$-quotient simply means that the multisets
\[
\rmst{\quo i\la}{i\in\calx_r}\qquad\text{and}\qquad\rmst{\quo i\mu}{i\in\calx_r}
\]
are equal for each $r$.

Now since $\xi$ is an $(s,t)$-core, it follows from the proof of \cite[Proposition 4.1]{mfcores} that the $\es i\xi$ lying in a given congruence class modulo $t$ form an arithmetic progression with common difference $t$. So given $r$, there are $a\in\zsz$ and $m\in\bbn_0$ such that
\[
\calx_r=a+t,a+2t,\dots,a+mt
\]
and there is an integer $c$ such that
\[
\es{a+bt}\xi=c+bt
\]
for $b=1,\dots,m$. Now given any $1<b\ls m$, let $i\in\zsz$ be such that $it-\ci st=a+bt$. Then (from the first paragraph of this proof) the effect of applying $w_i$ to $\la$ is to fix all the $\es i\la$, and to interchange $\quo{a+(b-1)t}\la$ and $\quo{a+bt}\la$, fixing all other parts of the $s$-quotient of $\la$. So by applying elements of $\weyl s$, we can re-order $\quo{a+t}\la,\dots,\quo{a+mt}\la$ arbitrarily without affecting the $s$-set of $\la$ or the rest of $\qou s\la$. By doing this for every $r$, we can apply elements of $\weyl s$ to transform $\la$ into~$\mu$.
\end{pf}

\section{Generalised cores}\label{gcsec}

Equipped with our definitions and basic results concerning the action of $\weyl s$, we now come to our main object of study.

\subsection{Olsson's theorem}

We begin by stating a theorem of Olsson which is the starting point for the work in this paper; this says that the $t$-core of an $s$-core is again an $s$-core, but we phrase this slightly differently.

\begin{thmc}{ol}{Theorem 1}\label{olssonnath}
Suppose $s,t$ are coprime positive integers and $\la\in\calp$. Then
\[
\wei s\la=0\quad\Longrightarrow\quad \wei s{(\cor t\la)}=0.
\]
\end{thmc}

Now we can give our first main result, which is a generalisation of \cref{olssonnath}.

\begin{thm}\label{ineq}
Suppose $s,t$ are coprime positive integers and $\la\in\calp$.  Then
\[
\wei s{(\cor t\la)}\ls\wei s\la.
\]
\end{thm}

\begin{pf}
We proceed by induction on $\wei t\la$, with the case where $\la$ is a $t$-core being trivial. Assuming $\la$ is not a $t$-core, we can find $b\in\be\la$ such that $b-t\notin\be\la$. We define a new partition $\nu$ by replacing $a$ with $a-t$ for every $a\in\be\la$ such that $a-t\notin\be\la$ and $a\equiv b\ppmod s$.  Then $\cor t\nu=\cor t\la$ and $\wei t\nu<\wei t\la$, so by induction it suffices to show that $\wei s\nu\ls\wei s\la$.

We use \cref{weightbeta} to compare $\wei s\la$ and $\wei s\nu$. Call a pair $(x,l)\in\bbz\times\bbn$ a \emph{weight pair} for $\la$ if $x\in\be\la\notni x-ls$. If $x\nequiv b,b-t\ppmod s$, then clearly $(x,l)$ is a weight pair for $\la$ if and only if it is a weight pair for $\nu$. Now suppose $x\equiv b\ppmod s$, and consider the two pairs $(x,l)$ and $(x-t,l)$. By considering the sixteen possibilities for the set $\be\la\cap\{x,x{-}ls,x{-}t,x{-}ls{-}t\}$, we can check that among the two pairs $(x,l)$ and $(x{-}t,l)$, there are at least as many weight pairs for $\la$ as there are for $\nu$; hence $\wei s\nu\ls\wei s\la$.
\end{pf}

With any inequality, it is natural to consider the situation where equality occurs. Hence we make the following definition.

\begin{defn}
Suppose $s,t$ are positive integers.  A partition $\la$ is an \emph{\gc st} if
\[
\wei s{(\cor t\la)}=\wei s\la.
\]
We write  $\cg st$ for the set of \gc sts.
\end{defn}

Trivially $\cg st$ includes all $t$-cores, and it follows from \cref{olssonnath} that $\cg st$ contains all $s$-cores.  However, in general $\cg st$ will include partitions which are neither $s$- nor $t$-cores; for example, $(4,1)$ is a \gc23.

Note that in the above definition we do not assume that $s$ and $t$ are coprime. However, for the rest of this paper we do make this assumption. Given this, it is easy to determine whether a partition is an \gc st from its beta-set.

\begin{propn}\label{betacondition}
Suppose $\la\in\calp$, $r\in\bbz$ and $s,t$ are coprime positive integers. Then $\la$ is an \gc st if and only if there do not exist integers $d,e,f$ such that:
\begin{itemize}
\item
$d\equiv e\ppmod s$;
\item
$d\equiv f\ppmod t$;
\item
$d,e+f-d\in\ber r\la$;
\item
$e,f\notin\ber r\la$.
\end{itemize}
\end{propn}

\begin{pf}
Since $\ber r\la$ is just a translation of $\be\la$, we may assume $r=0$. Say that $(d,e,f)$ is a \emph{bad triple} for $\la$ if $d,e,f$ satisfy the conditions in the \lcnamecref{betacondition}. First we suppose $(d,e,f)$ is bad, and show that $\la$ is not a $t$-core. Trivially, we must have either $d>f$ or $e+f-d>e$, either way, we find that there are $x,y\in\bbz$ such that $x>y$, $x\equiv y\ppmod t$ and $x\in\be\la\notni y$. So by \cref{basicbeta} $\la$ is not an $(x-y)$-core, and hence by \cref{scoremscore} $\la$ is not a $t$-core.

So (since every $t$-core is an \gc st) the \lcnamecref{betacondition} is true for $\la$ a $t$-core. Now we assume $\la$ is not a $t$-core, and choose $b\in\be\la$ such that $b-t\notin\be\la$.  As in the proof of \cref{ineq}, we define a new partition $\nu$ by replacing $a$ with $a-t$ for every $a\in\be\la$ such that $a-t\notin\be\la$ and $a\equiv b\ppmod s$. Now by induction it suffices to show that either:
\begin{itemize}
\item
$\wei s\nu=\wei s\la$, and there is a bad triple for $\nu$ if and only if there is a bad triple for $\la$; or
\item
$\wei s\nu<\wei s\la$, and there is a bad triple for $\la$.
\end{itemize}

Suppose first that there is no $a\equiv b\ppmod s$ for which $a\notin\be\la\ni a-t$. Then we have $\nu=w_i\la$, where $i\in\zsz$ is such that $it-\ci st=b+s\bbz$, and $w_i$ is the corresponding generator of $\weyl s$ acting at level $t$. So by \cref{orbitprops}(\ref{op2}) we have $\wei s\nu=\wei s\la$; and $(d,e,f)$ is a bad triple for $\la$ if and only if $(w_id,w_ie,w_if)$ is a bad triple for $\nu$.

Next suppose there is an $a\equiv b\ppmod t$ for which $a\notin\be\la\ni a-t$. Then $(b,a,b-t)$ is a bad triple for $\la$, and it remains to show that $\wei s\nu<\wei s\la$. As in the proof of \cref{ineq}, we consider weight pairs $(x,l)$. Taking $x=\max\{a,b\}$ and $l=\mfrac{|a-b|}t$, we find that exactly one of $(x,l)$ and $(x-t,l)$ is a weight pair for $\la$, while neither of them is a weight pair for $\nu$. Using the rest of the argument in the proof of \cref{ineq}, we have $\wei s\nu<\wei s\la$.
\end{pf}

\begin{cory}\label{stts}
If $s,t$ are coprime positive integers, then $\cg st=\cg ts$.
\end{cory}

\begin{pf}
The condition in \cref{betacondition} is symmetric in $s$ and $t$.
\end{pf}

This last result (which is very surprising given the definition of $\cg st$) suggests that the set $\cg st$ is worth studying. Our intuition is that $\cg st$, rather than $\cores s\cup\cores t$, is the `correct' counterpart to $\cores s\cap\cores t$ (just as one studies the sum of two subspaces of a vector space rather than their union).

We now go on to examine the structure of $\cg st$ with respect to the level $t$ action of $\weyl s$.

\subsection{Affine symmetric group actions on $\cg st$}

We continue to assume that $s$ and $t$ are coprime. Recall that the group $\weyl s$ acts at level $t$ on $\calp$; symmetrically, $\weyl t$ acts at level $s$ on $\calp$. These actions commute (since the actions on $\bbz$ commute), and so we have an action of $\weyl s\times\weyl t$ on $\calp$.

Our first result is that $\cg st$ is a union of orbits for the level $t$ action of $\weyl s$.

\begin{propn}\label{unorbs}
Suppose $s$ and $t$ are coprime positive integers. Given $\la\in\calp$ and $w\in\weyl s$, define $w\la$ using the level $t$ action. If $\la\in\cg st$, then $w\la\in\cg st$.
\end{propn}

\begin{pf}
Using \cref{orbitprops}(\ref{op1},\ref{op2}) and the fact that $\la\in\cg st$, we have
\[
\wei s{(\cor t{(w\la)})}=\wei s{(\cor t\la)}=\wei s\la=\wei s{(w\la)}.\qedhere
\]
\end{pf}

Interchanging $s$ and $t$ and appealing to \cref{stts}, we see that $\cg st$ is also a union of orbits for the level $s$ action of $\weyl t$. Hence $\cg st$ is a union of orbits for the action of $\weyl s\times\weyl t$.

Now we consider these orbits in more detail. We begin by considering just the level $t$ action of $\weyl s$.

\begin{propn}\label{coreinorbit}
Suppose $\la\in\calp$, and let $\calo$ be the orbit containing $\la$ under the level $t$ action of $\weyl s$.  Then the following are equivalent.
\begin{enumerate}
\item\label{lagc}
$\la$ is an \gc st.
\item\label{concor}
$\calo$ contains a $t$-core.
\item\label{concorla}
$\calo$ contains $\cor t\la$.
\end{enumerate}
\end{propn}

\begin{pf}
Since every $t$-core is an \gc st, \cref{unorbs} shows that if $\calo$ contains a $t$-core, then $\la$ is an \gc st. So (\ref{concor}) implies (\ref{lagc}). Trivially (\ref{concorla}) implies (\ref{concor}), so it remains to show that (\ref{lagc}) implies (\ref{concorla}).

So suppose $\la$ is an \gc st{}. We can assume that $\la$ is not a $t$-core, so there is $b\in\be\la$ such that $b-t\notin\be\la$. From the proof of \cref{betacondition}, there is no $a\equiv b\ppmod s$ for which $a-t\in\be\la\notni a$, and if we take $i\in\zsz$ such that $it-\ci st=b+s\bbz$, then the partition $\nu=w_i\la$ satisfies $\cor t\nu=\cor t\la$ and $\wei t\nu<\wei t\la$. By \cref{orbitprops} $\nu$ is also an \gc st, and by induction the orbit containing $\nu$ contains $\cor t\nu$.
\end{pf}

Now we can introduce a connection between \gc sts and $(s,t)$-cores.

\begin{cory}\label{exone}
Let $\calo$ be an orbit of $\weyl s\times\weyl t$ consisting of \gc sts.  Then $\calo$ contains exactly one $(s,t)$-core.
\end{cory}

\begin{pf}
Let $\la$ be a partition in $\calo$.  Then by \cref{coreinorbit} $\cor t\la\in\calo$, and by the same result with $s$ and $t$ interchanged, $\nu:=\cor s{(\cor t\la)}$ lies in $\calo$. By \cref{olssonnath} $\nu$ is an $(s,t)$-core.

Now suppose that there is another $(s,t)$-core in $\calo$.  We can write this as $xw\nu$, with $w\in\weyl s$ and $x\in\weyl t$.  By \cref{orbitprops}(\ref{op2}) we have
\[
\wei t{w\nu}=\wei t{xw\nu}=0,
\]
so
\[
w\nu=\cor t{w\nu}=\cor t\nu=\nu,
\]
using \cref{orbitprops}(\ref{op1}).  Similarly $x\nu=\nu$, and so $xw\nu=x\nu=\nu$.
\end{pf}

\begin{rmks}\indent
\begin{enumerate}
\vspace{-\topsep}
\item
From \cref{coreinorbit,exone} we see that two \gc sts $\la$ and $\mu$ lie in the same orbit of $\weyl s\times\weyl t$ if and only if $\cor s{(\cor t\la)}=\cor s{(\cor t)}\mu$. But it does not seem to be easy to tell when two arbitrary partitions lie in the same orbit; we would like an analogue of \cref{samewquot}, but the author has so far been unable to find one.
\item
\cref{exone} shows that the number of orbits of $\weyl s\times\weyl t$ consisting of \gc sts equals the number of $(s,t)$-cores; by Anderson's theorem \cite[Theorems 1 \& 3]{and} this is exactly $\mfrac1{s+t}\mbinom{s+t}s$.
\end{enumerate}
\end{rmks}

\subsection{Further properties of \gc sts}

Now we give two more properties of \gc sts which will be useful later.

\begin{lemma}\label{basicprops}
Suppose $\la$ is an \gc st.  Then:
\begin{enumerate}
\item
$\la$ is an $st$-core;
\item\label{bp2}
$\cor s{(\cor t\la)}=\cor t{(\cor s\la)}$.
\end{enumerate}
\end{lemma}

\begin{pfenum}
\item
If $\la$ is not an $st$-core, let $\nu$ be a partition obtained by removing an $st$-hook.  Then (from \cref{basicbeta}) we can obtain $\nu$ from $\la$ either by successively removing $s$ $t$-hooks, or by successively removing $t$ $s$-hooks.  Hence we see that
\[
\cor t\nu=\cor t\la,\qquad \wei s\nu=\wei s\la-t.
\]
So
\[
\wei s{(\cor t\la)}=\wei s{(\cor t\nu)}\ls\wei s\nu<\wei s\la,
\]
so $\la$ is not an \gc st.
\item
By \cref{olssonnath} both $\cor s{(\cor t\la)}$ and $\cor t{(\cor s\la)}$ are $(s,t)$-cores, and from \cref{coreinorbit} they both lie in the same orbit as $\la$ under the action of $\weyl s\times\weyl t$.  So the result follows from \cref{exone}.
\end{pfenum}

Note that the properties in \cref{basicprops} do not characterise \gc sts; for example, the rectangular partition $(s^t)$ satisfies both properties, but if $s,t>1$ it is not an \gc st.

\section{The sum of an $s$-core and a $t$-core}\label{sumsec}

Next we consider the possibility of constructing an \gc st{} with specified $s$-core and $t$-core.  We continue to assume that $s$ and $t$ are coprime.

\subsection{Partitions with a given $s$-core and $t$-core}

It is a simple exercise using the Chinese Remainder Theorem to show that given an $s$-core $\sigma$ and a $t$-core $\tau$, there are infinitely many partitions $\la$ with $\cor s\la=\sigma$ and $\cor t\la=\tau$.  But if we insist that $\la$ be an \gc st, then by \cref{basicprops}(\ref{bp2}) we need $\cor t\sigma=\cor s\tau$.  In this case, we have the following result.

\begin{propn}\label{speccore}
Suppose $\sigma\in\cores s$ and $\tau\in\cores t$, and that $\cor t\sigma=\cor s\tau$.  Then there is a unique \gc st{} $\la$ such that $\cor s\la=\sigma$ and $\cor t\la=\tau$.  Moreover, $|\la|=|\sigma|+|\tau|-|\cor s\tau|$, and $\la$ is the unique smallest partition with $s$-core $\sigma$ and $t$-core $\tau$.
\end{propn}

\begin{pf}
In this proof we use \cref{orbitprops} without comment.  Let $\xi=\cor t\sigma$, and consider the action of $\weyl s\times \weyl t$ on $\calp$.  By \cref{coreinorbit} we can find $w\in\weyl s$ and $x\in\weyl t$ such that $w\xi=\sigma$ and $x\xi=\tau$, and we let $\la=w\tau$.  Then $\la$ is an \gc st, since it lies in the same orbit as $\tau$. Moreover, $\cor t\la=\cor t\tau=\tau$, and
\[
\cor s\la=\cor s{(wx\xi)}=\cor s{(xw\xi)}=\cor s{(x\sigma)}=\cor s\sigma=\sigma.
\]
Furthermore,
\begin{align*}
|\la|&=|\cor s\la|+s\wei s\la\\
&=|\sigma|+s\wei s{(w\tau)}\\
&=|\sigma|+s\wei s{\tau}\\
&=|\sigma|+|\tau|-|\cor s\tau|.
\end{align*}

Now suppose $\mu$ is a partition other than $\la$ with $s$-core $\sigma$ and $t$-core $\tau$, and let $w,x$ be as above.  Then we have $\cor t{(w^{-1}\mu)}=\cor t\mu=\tau$, but $w^{-1}\mu\neq w^{-1}\la=\tau$.  So $|w^{-1}\mu|>|\tau|$.  Hence
\[
\wei s\mu=\wei s{(w^{-1}\mu)}=\frac{|w^{-1}\mu|-|\xi|}s>\frac{|\tau|-|\xi|}s=\wei s\tau;
\]
so $\mu$ is not an \gc st. Furthermore, we see that $\wei s\mu>\wei s \la$, so $|\mu|>|\la|$, and hence $|\la|$ is the unique smallest partition with $s$-core $\sigma$ and $t$-core $\tau$.
\end{pf}

We write $\plu\sigma\tau$ for the partition $\la$ given by \cref{speccore}.

\begin{rmk}
We have shown that $\plu\sigma\tau$ is the smallest partition with $s$-core $\sigma$ and $t$-core $\tau$ in terms of size; it is reasonable to ask whether $\plu\sigma\tau$ is smallest in the sense that $\plu\sigma\tau\subseteq\mu$ for any partition $\mu$ with $s$-core $\sigma$ and $t$-core $\tau$. In fact, this is false: taking $(s,t)=(2,3)$, we have $\plu{(2,1)}{(2)}=(2,1^3)$; but the partition $(8,3)\nsupseteq(2,1^3)$ also has $2$-core $(2,1)$ and $3$-core $(2)$.
\end{rmk}

Now we derive a corollary which yields another characterisation of \gc sts.

\begin{cory}\label{uniquest}
Suppose $s$ and $t$ are coprime positive integers not both equal to $1$ and $\la$ is a partition of $n$. Then $\la$ is an \gc st{} if and only if there is no other partition of $n$ with the same $s$-core and $t$-core as $\la$.
\end{cory}

\begin{pf}
Suppose $\la$ is an \gc st, and let $\sigma=\cor s\la$ and $\tau=\cor t\la$. Then $\la=\plu\sigma\tau$, since by \cref{speccore} this is the unique \gc st{} with $s$-core $\sigma$ and $t$-core $\tau$. So by the last part of \cref{speccore}, $\la$ is the unique partition of its size with $s$-core $\sigma$ and $t$-core $\tau$.

Conversely, suppose $\la$ is not an \gc st. Then by \cref{betacondition} we can find integers $d,e,f$ such that $d\equiv e\ppmod s$, $d\equiv f\ppmod t$ and $d,e+f-d$ are elements of $\be\la$ while $e,f$ are not.

First assume $d,e,f$ and $e+f-d$ are distinct. Define a new partition $\mu$ by
\[
\be\mu=\be\la\setminus\{d,e+f-d\}\cup\{e,f\}.
\]
Then by \cref{basicbeta} $\mu$ can obtained from $\la$ by removing a $|d-e|$-hook and adding a $|d-e|$-hook. So $|\mu|=|\la|$, and $\mu$ has the same $s$-core as $\la$, since $s$ divides $d-e$. Alternatively, $\mu$ can be obtained from $\la$ by removing a $|d-f|$-hook and adding a $|d-f|$-hook, so $\mu$ also has the same $t$-core as $\la$. So $\la$ is not the unique partition of $n$ with $s$-core $\sigma$ and $t$-core $\tau$.

Now assume $d,e,f$ and $e+f-d$ are not distinct. Then these integers are congruent modulo $st$, and hence by \cref{basicbeta} $\la$ is not an $st$-core. Since $st>1$ by assumption, it is easy to find another partition $\mu$ of $n$ with the same $st$-core as $\la$, and by \cref{scoremscore} $\mu$ has the same $s$-core and $t$-core as $\la$.
\end{pf}

\subsection{Constructing $\plu\sigma\tau$}

Suppose $\sigma$ is an $s$-core and $\tau$ is a $t$-core, with $\cor t\sigma=\cor s\tau$. In this section we give a method for constructing the partition $\plu\sigma\tau$. Of course, this can be done as in the proof of \cref{speccore}: find $w\in\weyl s$ such that $\sigma=w\cor t\sigma$, and then compute $w\tau$. But this is a laborious process; we present here a much quicker method using weighted quotients.

Recall the $s$-set $\lset{\es i\la}{i\in\zsz}$ and the $s$-quotient $\rtup{\quo i\la}{i\in\zsz}$ of a partition $\la$.

\begin{lemma}\label{samequot}
Suppose $j,k\in\zsz$ and $\tau$ is a $t$-core with $\es j\tau\equiv\es k\tau\ppmod t$. Then $\quo j\tau=\quo k\tau$.
\end{lemma}

\begin{pf}
Without loss of generality suppose $\es k\tau>\es j\tau$. The elements of $\be\tau\cap j$ are the integers $s(\quo j\tau_i-i)+\es j\tau$, for $i\gs1$, and similarly for $\be\tau\cap k$. If $\quo j\tau\neq\quo k\tau$, then by \cref{becont} $\be{\quo k\tau}\nsubseteq\be{\quo j\tau}$, so there is some $i$ such that $\quo k\tau_i-i\notin\be{\quo j\tau}$. This means that $s(\quo k\tau_i-i)+\es j\tau\notin\be\tau$; on the other hand, $(\quo k\tau_i-i)s+\es k\tau\in\be\tau$, so by \cref{basicbeta} $\tau$ is not an $(\es k\tau-\es j\tau)$-core, and hence by \cref{scoremscore} $\tau$ is not a $t$-core. Contradiction.
\end{pf}

Now recall the $t$-weighted $s$-quotient $\rmst{(\es i\la+t\bbz,\quo i\la)}{i\in\zsz}$ of $\la$.

\begin{propn}\label{constructing}
Suppose $\sigma$ is an $s$-core and $\tau$ is a $t$-core and that $\cor t\sigma=\cor s\tau$. Then there is a unique partition with $s$-core $\sigma$ and with the same $t$-weighted $s$-quotient as $\tau$.
\end{propn}

\begin{pf}
Let $\xi=\cor t\sigma$. Then by \cref{samequot,coreinorbit} $\xi$ has the same $t$-weighted $s$-quotient as $\sigma$, i.e.\ there is a bijection $\phi:\zsz\to\zsz$ such that $\es i\sigma\equiv \es{\phi(i)}\xi\ppmod t$ for each $i$. Since $\xi=\cor s\tau$, it is therefore possible to construct a partition $\la$ as stated: we just take the partition $\la$ with $s$-core $\sigma$, and with $\quo i\la=\quo {\phi(i)}\tau$ for each $i$.

By \cref{samequot} we have $\quo j\tau=\quo k\tau$ whenever $\es j\tau\equiv\es k\tau\ppmod t$, i.e.\ whenever $\es{\phi^{-1}(j)}\sigma\equiv\es{\phi^{-1}(k)}\sigma\ppmod t$, so we have no choice in the construction of $\la$, and $\la$ is unique.
\end{pf}

\begin{propn}\label{constructed}
Suppose $\sigma$ is an $s$-core and $\tau$ is a $t$-core, with $\cor t\sigma=\cor s\tau$, and let $\la$ be the partition with $s$-core $\sigma$ and with the same $t$-weighted $s$-quotient as $\tau$. Then $\la=\plu\sigma\tau$.
\end{propn}

\begin{pf}
By \cref{samewquot}, $\la$ and $\tau$ lie in the same level $t$ orbit of $\weyl s$. Hence by \cref{coreinorbit} $\la$ is an \gc st{} and $\cor t\la=\tau$. Since $\cor s\la=\sigma$ by construction, we have $\la=\plu\sigma\tau$.
\end{pf}

\begin{eg}
Take $s=3$, $t=5$, $\sigma=(7,5,4^2,3^2,2^2,1^2)$ and $\tau=(5,2^4,1^4)$. Then $\cor5\sigma=\cor3\tau=(2)$. We have
\[
(\es{3\bbz}\tau,\es{1+3\bbz}\tau,\es{2+3\bbz}\tau)=(0,4,-1),\qquad(\quo{3\bbz}\tau,\quo{1+3\bbz}\tau,\quo{2+3\bbz}\tau)=((1^3),(1),(1)),
\]
so that the $5$-weighted $3$-quotient of $\tau$ is
\[
[(5\bbz,(1^3)),(4+5\bbz,(1)),(4+5\bbz,(1))].
\]
On the other hand,
\[
(\es{3\bbz}\sigma,\es{1+3\bbz}\sigma,\es{2+3\bbz}\sigma)=(9,4,-10);
\]
so if $\la$ has $3$-core $\sigma$ (and hence has $s$-set $\{9,4,-10\}$), then the only way $\la$ can have the same $5$-weighted $3$-quotient as $\tau$ is if $(\quo{3\bbz}\la,\quo{1+3\bbz}\la,\quo{2+3\bbz}\la)=((1),(1),(1^3))$. This gives $\la=(10,6^2,4,3^2,2^2,1^{11})$. The $3$-runner abacus displays of $\sigma,\tau$ and $\la$ are as follows.
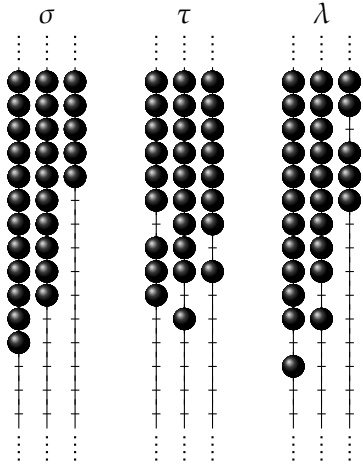
\[
\begin{array}{c@{\qquad}c@{\qquad}c}
\sigma&\tau&\la\\
\abacus(vvv,bbb,bbb,bbb,bbb,bbb,bbn,bbn,bbn,bbn,bbn,bnn,bnn,nnn,nnn,nnn;vvv)&
\abacus(vvv,bbb,bbb,bbb,bbb,bbb,bbb,nbb,bbn,bbb,bnn,nbn,nnn,nnn,nnn,nnn;vvv)&
\abacus(vvv,bbb,bbb,bbn,bbb,bbb,bbb,bbn,bbn,bbn,bnn,bbn,nnn,bnn,nnn,nnn;vvv)
\end{array}
\]
For further illustration, we take the same example with $s$ and $t$ interchanged. The $3$-weighted $5$-quotient of $\sigma$ is
\[
[(3\bbz,(1^2)),(3\bbz,(1^2)),(2+3\bbz,\varnothing),(3\bbz,(1^2)),(2+3\bbz,\varnothing)].
\]
On the other hand,
\[
(\es{5\bbz}\tau,\es{1+5\bbz}\tau,\es{2+5\bbz}\tau,\es{3+5\bbz}\tau,\es{4+5\bbz}\tau)=(5,-9,2,3,9);
\]
so we must take $(\quo{5\bbz}\la,\quo{1+5\bbz}\la,\quo{2+5\bbz}\la,\quo{3+5\bbz}\la,\quo{4+5\bbz}\la)=(\varnothing,(1^2),\varnothing,(1^2),(1^2))$, again yielding $\la=(10,6^2,4,3^2,2^2,1^{11})$. The $5$-runner abacus displays are as follows.
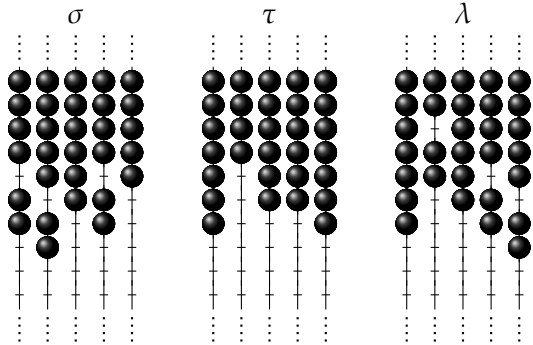
\[
\begin{array}{c@{\qquad}c@{\qquad}c}
\sigma&\tau&\la\\
\abacus(vvvvv,bbbbb,bbbbb,bbbbb,bbbbb,nbbnb,bnbbn,bbnbn,nbnnn,nnnnn,nnnnn;vvvvv)&
\abacus(vvvvv,bbbbb,bbbbb,bbbbb,bbbbb,bnbbb,bnbbb,bnnnb,nnnnn,nnnnn,nnnnn;vvvvv)&
\abacus(vvvvv,bbbbb,bbbbb,bnbbb,bbbbb,bbbnb,bnbbn,bnnbb,nnnnb,nnnnn,nnnnn;vvvvv)
\end{array}
\]
\end{eg}

\section{The $\kappa$-orbit}\label{kappasec}

In this section we examine one particular orbit in $\cg st$, continuing to assume that $s$ and $t$ are coprime. Under this assumption, there are only finitely many $(s,t)$-cores, and there is a unique largest such. This partition (which is usually denoted $\kappa_{s,t}$) has been studied before; it is known \cite[Theorem 4.1]{os} that $|\kappa_{s,t}|=\frac1{24}(s^2-1)(t^2-1)$, and also that if $\la$ is any $(s,t)$-core then $\la\subseteq\kappa_{s,t}$ \cite[Theorem 2.4]{vand}, \cite[Theorem 5.1]{mfcores}. In this section we will consider the $\weyl s\times\weyl t$-orbit containing $\kappa_{s,t}$. We denote this orbit $\cg st^\kappa$, and refer to it as the \emph{$\kappa$-orbit} of $\cg st$. We will see that $\cg st^\kappa$ is naturally in bijection with $\cores s\times\cores t$.

To begin with, we explain how to construct $\kappa_{s,t}$. Let $\calb_{s,t}$ denote the set of integers which cannot be written as a linear combination of $s$ and $t$ with non-negative integer coefficients; this set can be written as
\[
\calb_{s,t}=\lset{at-bs}{a\in\{0,\dots,s-1\},\ b\in\bbn}.
\]
Then $\calb_{s,t}$ is bounded above and its complement in $\bbz$ is bounded below, so it is a beta-set of a partition, and this partition is $\kappa_{s,t}$. In fact (recalling the integer $\ci st=\frac12(s-1)(t-1)$ from \cref{groupsec}) $\calb_{s,t}=\ber{\ci st}{\kappa_{s,t}}$.

The following statement is proved in \cite[\S5]{mfcores}.

\begin{lemma}\label{kappasset}
The $s$-set of $\kappa_{s,t}$ is
\[
\{-\ci st,t-\ci st,2t-\ci st,\dots,(s-1)t-\ci st\}.
\]
\end{lemma}

Note in particular that the elements of the $s$-set of $\kappa_{s,t}$ are congruent modulo $t$; in fact, $\kappa_{s,t}$ is the unique $(s,t)$-core with this property.

\begin{eg}
Take $s=3$ and $t=4$.  Then
\begin{align*}
\calb_{s,t}&=\{-3,-6,-9,\dots\}\cup\{1,-2,-5,\dots\}\cup\{5,2,-1,\dots\}\\
&=\{5,2,1,-1,-2,-3,\dots\}\\
&=\ber3{(3,1^2)},
\end{align*}
so $\kappa_{s,t}=(3,1^2)$. The $3$-set of this partition is $\{-3,1,5\}$, while its $4$-set is $\{-3,0,3,6\}$.
\end{eg}

Now we consider the $\kappa$-orbit $\cg st^\kappa$. We begin by showing that the $s$- and $t$-quotients of partitions in this orbit have a particularly nice form.

\begin{propn}\label{homogcore}
Suppose $\tau$ is a $t$-core such that $\cor s\tau=\kappa_{s,t}$. Then $\qou s\tau=(\la,\dots,\la)$ for some $t$-core~$\la$.
\end{propn}

\begin{pf}
By \cref{kappasset}, the elements of the $s$-set of $\tau$ (i.e.\ the $s$-set of $\kappa_{s,t}$) are congruent modulo $t$. By \cref{samequot}, this means that the components of $\qou s\tau$ are all equal. Furthermore, since $\tau$ is a $t$-core, it is an $st$-core, by \cref{scoremscore}, and so by \cref{quotprops}(\ref{wc1}) each component of $\qou s\tau$ must be a $t$-core.
\end{pf}

\begin{rmk}
Let us say that a partition is \emph{$s$-homogeneous} if all the components of its $s$-core are equal; we have just shown that a $t$-core whose $s$-core is $\kappa_{s,t}$ is $s$-homogeneous. However, the condition $\cor s\tau=\kappa_{s,t}$ is not necessary for a $t$-core $\tau$ to be $s$-homogeneous; for example, $\tau=\varnothing$ has $s$-quotient $(\varnothing,\dots,\varnothing)$. However, one can show that if $t$ is prime, then there are only finitely many $s$-homogeneous $t$-cores whose $s$-core is not $\kappa_{s,t}$.
\end{rmk}

A consequence of \cref{homogcore2} is that the construction of $\plu\sigma\tau$ is even simpler when $\sigma$ is an $s$-core and $\tau$ a $t$-core in the $\kappa$-orbit.

\begin{propn}\label{stsumkappa}
Suppose $\sigma$ is an $s$-core and $\tau$ a $t$-core with $\cor t\sigma=\cor s\tau=\kappa_{s,t}$. Then $\plu\sigma\tau$ is the partition with $s$-core $\sigma$ and the same $s$-quotient as $\tau$.
\end{propn}

\begin{pf}
$\plu\sigma\tau$ has $s$-core $\sigma$ by definition.  By \cref{coreinorbit} $\tau$ and $\plu\sigma\tau$ lie in the same level $t$ orbit of $\weyl s$, and so by \cref{orbitprops}(\ref{op1a}) have the same $s$-quotient up to re-ordering.  Since $\tau$ has $s$-quotient $(\la,\dots,\la)$ for some $\la$, $\plu\sigma\tau$ does too. 
\end{pf}

\begin{eg}
Take $(s,t)=(3,4)$, so that $\kappa_{s,t}=(3,1^2)$. Then $\sigma=(15,13,11,9,7,5,3,2^2,1^2)$ is a $3$-core with $4$-core $(3,1^2)$, while $\tau=(6,4^2,2^3,1^3)$ is a $4$-core with $3$-core $(3,1^2)$. Hence $\plu\sigma\tau$ is the partition with $3$-core $\sigma$ and the same $3$-quotient as $\tau$, or equivalently the partition with $4$-core $\tau$ and the same $4$-quotient as $\sigma$. We have
\[
\qou4\sigma=((3,1),(3,1),(3,1),(3,1)),\qquad\qou3\tau=((1^2),(1^2),(1^2)),
\]
and by either route we find that $\plu\sigma\tau=(18,16,11,9,7,5,4^2,3^2,1^7)$. The $3$- and $4$-runner abacus displays for these partitions are as follows.
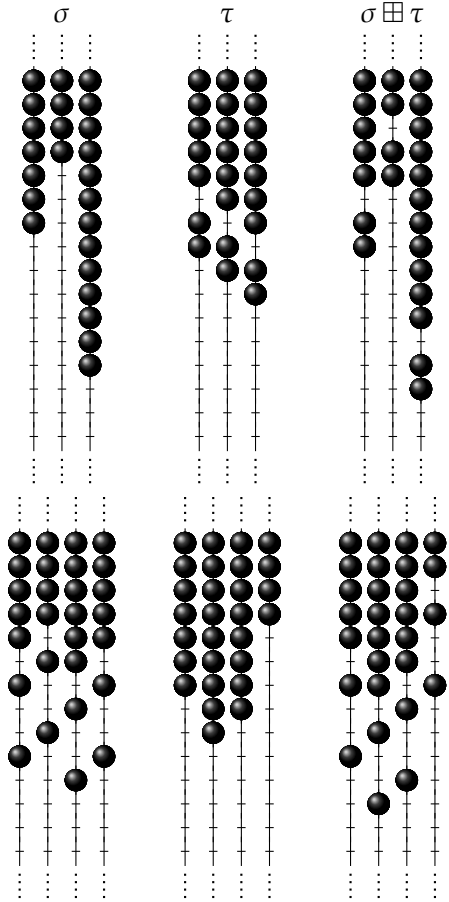
\[
\begin{array}{c@{\qquad}c@{\qquad}c}
\sigma&\tau&\plu\sigma\tau\\
\abacus(vvv,bbb,bbb,bbb,bbb,bnb,bnb,bnb,nnb,nnb,nnb,nnb,nnb,nnb,nnn,nnn,nnn;vvv)&
\abacus(vvv,bbb,bbb,bbb,bbb,bbb,nbb,bnb,bbn,nbb,nnb,nnn,nnn,nnn,nnn,nnn,nnn;vvv)&
\abacus(vvv,bbb,bbb,bnb,bbb,bbb,nnb,bnb,bnb,nnb,nnb,nnb,nnn,nnb,nnb,nnn,nnn;vvv)\\
\abacus(vvvv,bbbb,bbbb,bbbb,bbbb,bnbb,nbbn,bnnb,nnbn,nbnn,bnnb,nnbn,nnnn,nnnn,nnnn;vvvv)&
\abacus(vvvv,bbbb,bbbb,bbbb,bbbb,bbbn,bbbn,bbbn,nbbn,nbnn,nnnn,nnnn,nnnn,nnnn,nnnn;vvvv)&
\abacus(vvvv,bbbb,bbbb,bbbn,bbbb,bbbn,nbbn,bbnb,nnbn,nbnn,bnnn,nnbn,nbnn,nnnn,nnnn;vvvv)
\end{array}
\]
\end{eg}

Now we give the converse to \cref{homogcore}.

\begin{propn}\label{homogcore2}
Suppose $\upsilon$ is a $t$-core, and let $\tau$ be the partition with $s$-core $\kappa_{s,t}$ and $s$-quotient $(\upsilon,\dots,\upsilon)$. Then $\tau$ is a $t$-core.
\end{propn}

\begin{pf}
We consider the beta-set $\ber{\ci st}\tau$. By \cref{basicbeta}, we must show that $c-t\in\ber{\ci st}\tau$ for every $c\in\ber{\ci st}\tau$. In other words, we must show that for each $b\in\zsz$ we have
\[
\lset{c-t}{c\in\ber{\ci st}\tau\cap b}\subseteq\ber{\ci st}\tau\cap(b-t).
\]
Since $s$ and $t$ are coprime, we can write $b=at+s\bbz$ for some $a\in\{0,\dots,s-1\}$. Then by construction we have
\[
\ber{\ci st}\tau\cap b=\lset{at+s(\upsilon_i-i)}{i\in\bbn}.
\]
If $a>0$, then we also have
\[
\ber{\ci st}\tau\cap(b-t)=\lset{(a-1)t+s(\upsilon_i-i)}{i\in\bbn},
\]
and therefore
\[
\lset{c-t}{c\in\ber{\ci st}\tau\cap b}=\ber{\ci st}\tau\cap(b-t).
\]
So it remains to consider the case $a=0$. Now we have
\[
\ber{\ci st}\tau\cap b=\lset{s(\upsilon_i-i)}{i\in\bbn},\qquad\ber{\ci st}\tau\cap(b-t)=\lset{(s-1)t+s(\upsilon_i-i)}{i\in\bbn}.
\]
If the left-hand side is not contained in the right-hand side, then for some $j$ we have
\[
s(\upsilon_j-j)-t\notin\lset{(s-1)t+s(\upsilon_i-i)}{i\in\bbn},
\]
which gives $\upsilon_j-j-t\notin\lset{\upsilon_i-i}{i\in\bbn}$. But by \cref{basicbeta}, this contradicts the fact that $\upsilon$ is a $t$-core.
\end{pf}

Now we can give a concrete description of $\cg st^\kappa$, which says that as a $\weyl s\times\weyl t$-set, $\cg st^\kappa$ is isomorphic to the product of $\cores s$ and $\cores t$ (with $\weyl s$ and $\weyl t$ acting at level $1$ on these factors).

\begin{propn}\label{bijec}\indent
\begin{enumerate}
\vspace{-\topsep}
\item
There is a bijection
\begin{align*}
\cg st^\kappa&\longrightarrow\cores s\times\cores t\\
\la&\longmapsto(\quo{t\bbz}\la,\quo{s\bbz}\la).
\end{align*}
\item
Given $\la\in\cg st^\kappa$ and $w\in\weyl s$, we have
\[
\quo{t\bbz}{(w\la)}=w(\quo{t\bbz}\la),
\]
where $w$ acts at level $t$ on $\cg st^\kappa$, and at level $1$ on $\cores s$.
\end{enumerate}
\end{propn}

\begin{pfenum}
\item
Let $\Theta$ denote the given map, and suppose $\la\in\cg st^\kappa$. Since $\la$ is an $st$-core, every component of its $s$-quotient is a $t$-core, and every component of its $t$-quotient is an $s$-core, by \cref{scoremscore}; so $\Theta$ really does map to $\cores s\times\cores t$. To show that $\Theta$ is a bijection, we construct an inverse. If $(\rho,\upsilon)\in\cores s\times\cores t$, let $\tau$ be the partition with $s$-core $\kappa_{s,t}$ and $s$-quotient $(\upsilon,\dots,\upsilon)$. Then $\tau$ is a $t$-core by \cref{homogcore2}. Similarly the partition $\sigma$ with $t$-core $\kappa_{s,t}$ and $t$-quotient $(\rho,\dots,\rho)$ is an $s$-core, and $\cor t\sigma=\cor s\tau$, so we can define a partition $\plu\sigma\tau$, which will lie in $\cg st^\kappa$. We define $\Xi(\rho,\upsilon)=\plu\sigma\tau$, and we have a function $\Xi:\cores s\times\cores t\to \cg st^\kappa$.

Now we show that $\Theta$ and $\Xi$ are mutual inverses. Suppose $\la\in\cg st^\kappa$, write $\Theta(\la)=(\rho,\upsilon)$ and let $\tau=\cor t\la$. Then $\tau$ is a $t$-core with $s$-core $\kappa_{s,t}$, so by \cref{homogcore} all the components of the $s$-quotient of $\tau$ are equal. But $\la$ and $\tau$ lie in the same level $t$ orbit of $\weyl s$, so have the same $s$-quotient up to re-ordering; since $\upsilon=\quo{s\bbz}\la$, this means that both $\la$ and $\tau$ have $s$-quotient $(\upsilon,\dots,\upsilon)$. So $\tau$ is the (unique) partition with $s$-core $\kappa_{s,t}$ and $s$-quotient $(\upsilon,\dots,\upsilon)$. Similarly the $s$-core $\sigma$ of $\la$ is the partition with $t$-core $\kappa_{s,t}$ and $t$-quotient $(\rho,\dots,\rho)$. So $\Xi(\Theta(\la))=\Xi(\rho,\upsilon)=\plu\sigma\tau$, which is the unique partition in $\cg st$ with $s$-core $\sigma$ and $t$-core $\tau$, i.e.\ $\la$.

Now take $(\rho,\upsilon)\in\cores s\times\cores t$. Let $\tau$ be the partition with $s$-core $\kappa_{s,t}$ and $s$-quotient $(\upsilon,\dots,\upsilon)$, and $\sigma$ the partition with $t$-core $\kappa_{s,t}$ and $t$-quotient $(\rho,\dots,\rho)$, so that $\Xi(\rho,\upsilon)$ is by definition $\plu\sigma\tau$. This partition has the same $t$-quotient as $\sigma$ since it lies in the same level $s$ orbit of $\weyl t$, and in particular each component of its $t$-quotient is $\rho$. Similarly each component of the $s$-quotient of $\plu\sigma\tau$ is $\upsilon$, and so $\Theta(\Xi(\rho,\upsilon))=\Theta(\plu\sigma\tau)=(\rho,\upsilon)$.
\item
It suffices to consider the case where $w=w_i$, for $i\in\zsz$. Write $j=it-\ci st$, and let $\sigma=\quo{t\bbz}\la$. Then by the definition of $t$-quotient,
\[
\be\la\cap t\bbz=\lset{bt+\es{t\bbz}\la}{b\in\be\sigma}.
\]
Now the definition of the level $t$ action of $\weyl s$ gives
\[
\be{w_i\la}=\lset{c+t}{c\in\be\la\cap(j-t)}\cup\lset{c-t}{c\in\be\la\cap j}\cup(\be\la\setminus(j-t\cup j)).
\]
From \cref{kappasset} (with $s$ and $t$ interchanged) we have $\es{t\bbz}\la\equiv-\ci st\ppmod s$. Hence for $b\in\bbz$ we have $bt+\es{t\bbz}\la\in j-t$ if and only if $b\in i-1$, while $bt+\es{t\bbz}\la\in j$ if and only if $b\in i$. So
\[
\be{w_i\la}=\rset{(w_ib)t+\es{t\bbz}\la}{b\in\be\sigma}=\rset{bt+\es{t\bbz}\la}{b\in\be{w_i\sigma}},
\]
where $\weyl s$ acts at level $t$ in the first term, and at level $1$ in the other two terms. Hence
\[
\quo{t\bbz}{(w_i\la)}=w_i\sigma,
\]
as required.
\end{pfenum}

Of course, part (2) of the \lcnamecref{bijec} also holds with $s$ and $t$ interchanged, yielding the desired statement about the action of $\weyl s\times\weyl t$.

\begin{figure}[p]\label{23fig}
\newcommand\drawbw[1]{\draw[line width=4pt,white]{#1};\draw{#1};}
\newcommand\drawbwd[1]{\draw[line width=4pt,white,dashed]{#1};\draw[dashed]{#1};}
{\footnotesize\[\tdplotsetmaincoords{65}{18}
\begin{tikzpicture}[tdplot_main_coords,scale=0.05]
\tikzstyle{every node}=[fill=white,inner sep=0.5pt]
\coordinate(bb)at(1cm,1.73cm,0cm);
\coordinate(cc)at(1cm,-1.73cm,0cm);
\coordinate(aa)at($(bb)+(cc)$);
\coordinate(dd)at(0cm,0cm,-3.8cm);
\coordinate(sb)at($.5*(bb)$);
\coordinate(sc)at($.5*(cc)$);
\coordinate(sa)at($.5*(aa)$);
\coordinate(sd)at($.5*(dd)$);
\draw($-1.2*(aa)$)node{$\varnothing$}--++(sd)node[diamond,draw,thin]{\scriptsize0}--++(sd)node{$(1)$}--++(sd)node[diamond,draw,thin]{\scriptsize1}--++(sd)node{$(2,1)$}--++($.25*(dd)$);
\draw[dashed]($-1.2*(aa)+2.25*(dd)$)--++($.25*(dd)$);
\foreach\x in {-1}{
\drawbw{($\x*(dd)$)--++(sa)node[circle,draw,thin,inner sep=.7pt]{\scriptsize0}--++(sa)--++(sb)node[circle,draw,thin,inner sep=.7pt]{\scriptsize1}--++(sb)--++(sa)node[circle,draw,thin,inner sep=.7pt]{\scriptsize2}--++(sa)--++(sb)node[circle,draw,thin,inner sep=.7pt]{\scriptsize0}--++(sb)--++($.25*(aa)$)}
\drawbwd{($\x*(dd)+2.25*(aa)+2*(bb)$)--++($.25*(aa)$)}
\draw($\x*(dd)$)node{$\varnothing$}++(aa)++(bb)node{$(2)$}++(aa)++(bb)node{$(4,2)$};
\drawbw{($\x*(dd)+(aa)$)--++(sc)node[circle,draw,thin,inner sep=.7pt]{\scriptsize2}--++(sc)}
\drawbw{($\x*(dd)+(aa)+(aa)+(bb)$)--++(sc)node[circle,draw,thin,inner sep=.7pt]{\scriptsize1}--++(sc)}
\draw($\x*(dd)+(aa)$)node{$(1)$}++(aa)++(bb)node{$(3,1)$};
\drawbw{($\x*(dd)+(aa)+(cc)$)--++(sa)node[circle,draw,thin,inner sep=.7pt]{\scriptsize1}--++(sa)--++(sb)node[circle,draw,thin,inner sep=.7pt]{\scriptsize2}--++(sb)--++($.25*(aa)$)}
\drawbwd{($\x*(dd)+2.25*(aa)+(cc)+(bb)$)--++($.25*(aa)$)}
\draw($\x*(dd)+(aa)+(cc)$)node{$(1^2)$}++(aa)++(bb)node{$(3,1^2)$};
\drawbw{($\x*(dd)+2*(aa)+(cc)$)--++(sc)node[circle,draw,thin,inner sep=.7pt]{\scriptsize0}--++(sc)}
\draw($\x*(dd)+2*(aa)+(cc)$)node{$(2,1^2)$};
\drawbw{($\x*(dd)+2*(aa)+2*(cc)$)--++($.25*(aa)$)}
\drawbwd{($\x*(dd)+2.25*(aa)+2*(cc)$)--++($.25*(aa)$)}
\draw($\x*(dd)+2*(aa)+2*(cc)$)node{$(2^2,1^2)$};
}
\drawbw{(0,0,0)--++(sa)node[circle,draw,thin,inner sep=.7pt]{\scriptsize0}--++(sa)--++(sb)node[circle,draw,thin,inner sep=.7pt]{\scriptsize1}--++(sb)--++(sa)node[circle,draw,thin,inner sep=.7pt]{\scriptsize2}--++(sa)--++(sb)node[circle,draw,thin,inner sep=.7pt]{\scriptsize0}--++(sb)--++($.25*(aa)$)}
\foreach\x in{1,2}{
\drawbw{($\x*(dd)$)--++(aa)--++(bb)--++(aa)--++(bb)--++($.25*(aa)$)}
}
\foreach\x in{0,1,2}{
\drawbwd{($\x*(dd)+2.25*(aa)+2*(bb)$)--++($.25*(aa)$)}
}
\drawbw{(0,0,0)--++(sd)node[diamond,draw,thin]{\scriptsize0}--++(dd)node[diamond,draw,thin]{\scriptsize1}--++($.75*(dd)$)}
\foreach\x in{1,2}{
\drawbw{($\x*(aa)+\x*(bb)$)--++($2.25*(dd)$)}
}
\foreach\x in{0,1,2}{
\drawbwd{($\x*(aa)+\x*(bb)$)++($2.25*(dd)$)--++($.25*(dd)$)}
}
\foreach\x in{0,1}{
\drawbw{($(aa)+\x*(aa)+\x*(bb)$)--++($2.25*(dd)$)}
\drawbwd{($(aa)+\x*(aa)+\x*(bb)$)++($2.25*(dd)$)--++($.25*(dd)$)}
}
\draw(0,0,0)node{$(1)$}++(aa)++(bb)node{$(5,3,1)$}++(aa)++(bb)node{$(9,7,5,3,1)$};
\draw(dd)node{$(4,3,2,1)$}++(aa)++(bb)node{$(8,3,2^3,1)$}++(aa)++(bb)node{$(12,7,5,3,2^3,1)$};
\draw(dd)++(dd)node{$(7,6,5,4,3,2,1)$}++(aa)++(bb)node{$(11,6,5,4,3,2^3,1)$}++(aa)++(bb)node{$(15,10,5,4^3,3,2^3,1)$};
\drawbw{(aa)--++(sc)node[circle,draw,thin,inner sep=.7pt]{\scriptsize2}--++(sc)}
\drawbw{($2*(aa)+(bb)$)--++(sc)node[circle,draw,thin,inner sep=.7pt]{\scriptsize1}--++(sc)}
\foreach\x in{0,1}{
\foreach\y in{1,2}{
\drawbw{($\y*(dd)+(aa)+\x*(aa)+\x*(bb)$)--++(cc)}
}
}
\draw(aa)node{$(3,1^2)$}++(aa)++(bb)node{$(7,5,3,1^2)$};
\draw(aa)++(dd)node{$(6,3,2,1^3)$}++(aa)++(bb)node{$(10,5,3^2,2,1^3)$};
\draw(aa)++(dd)++(dd)node{$(9,6,5,4,3,2,1^3)$}++(aa)++(bb)node{$(13,8,5,4,3^3,2,1^3)$};
\drawbw{($(aa)+(cc)$)--++(sa)node[circle,draw,thin,inner sep=.7pt]{\scriptsize1}--++(sa)--++(sb)node[circle,draw,thin,inner sep=.7pt]{\scriptsize2}--++(sb)--++($.25*(aa)$)}
\foreach\x in{1,2}{
\drawbw{($\x*(dd)+(aa)+(cc)$)--++(aa)--++(bb)--++($.25*(aa)$)}
}
\foreach\x in{0,1,2}{
\drawbwd{($\x*(dd)+2.25*(aa)+(cc)+(bb)$)--++($.25*(aa)$)}
}
\foreach\x in{0,1}{
\drawbw{($(aa)+(cc)+\x*(aa)+\x*(bb)$)--++($2.25*(dd)$)}
\drawbwd{($(aa)+(cc)+\x*(aa)+\x*(bb)+2.25*(dd)$)--++($.25*(dd)$)}
}
\draw(aa)++(cc)node{$(3,2^2,1^2)$}++(aa)++(bb)node{$(7,5,3,2^2,1^2)$};
\draw(dd)++(aa)++(cc)node{$(6,5,2,1^5)$}++(aa)++(bb)node{$(10,5,4^2,2,1^5)$};
\draw(dd)++(dd)++(aa)++(cc)node{$(9,8,5,4,3,2,1^5)$}++(aa)++(bb)node{$(13,8,7,4,3^3,2,1^5)$};
\drawbw{(aa)++(aa)++(cc)--++($2.25*(dd)$)}
\drawbwd{(aa)++(aa)++(cc)++($2.25*(dd)$)--++($.25*(dd)$)}
\drawbw{($2*(aa)+(cc)$)--++(sc)node[circle,draw,thin,inner sep=.7pt]{\scriptsize0}--++(sc)}
\foreach\x in{1,2}{
\drawbw{($\x*(dd)+2*(aa)+(cc)$)--++(cc)}
}
\draw(aa)++(aa)++(cc)node{$(5,3^2,2^2,1^2)$};
\draw(dd)++(aa)++(aa)++(cc)node{$(8,5,4,2^2,1^5)$};
\draw(dd)++(dd)++(aa)++(aa)++(cc)node{$(11,8,7,4,3,2^3,1^5)$};
\drawbw{(aa)++(aa)++(cc)++(cc)--++($2.25*(dd)$)}
\drawbwd{(aa)++(aa)++(cc)++(cc)++($2.25*(dd)$)--++($.25*(dd)$)}
\foreach\x in{0,1,2}{
\drawbw{($\x*(dd)+2*(aa)+2*(cc)$)--++($.25*(aa)$)}
\drawbwd{($\x*(dd)+2.25*(aa)+2*(cc)$)--++($.25*(aa)$)}
}
\draw(aa)++(aa)++(cc)++(cc)node{$(5,4^2,3^2,2^2,1^2)$};
\draw(dd)++(aa)++(aa)++(cc)++(cc)node{$(8,7,4,3^2,2^2,1^5)$};
\draw(dd)++(dd)++(aa)++(aa)++(cc)++(cc)node{$(11,10,7,6,3,2^5,1^5)$};
\end{tikzpicture}
\]}
\caption{The bijection between $\cores2\times\cores3$ and $\cg 23^\kappa$}
\end{figure}

In Figure 1, we illustrate \cref{bijec} in the case $(s,t)=(2,3)$. At the top of the diagram, we have drawn a portion of $\cores3$ as a labelled graph, with edges indicating the actions of the generators $w_{3\bbz},w_{1+3\bbz},w_{2+3\bbz}$ in the level $1$ action of $\weyl 3$. On the left, we have drawn a portion of $\cores2$, with edges representing the actions of $w_{2\bbz},w_{1+2\bbz}$ in the level $1$ action of $\weyl2$. The main part of the diagram shows a portion of $\cg23^\kappa$, which (we hope) makes the bijection $\cores2\times\cores3\to\cg23^\kappa$ clear. Here the edges represent the actions of the generators of $\weyl 3$ in the level $2$ action, and of $\weyl2$ in the level $3$ action.

\section{Final remarks}\label{rmksec}

\subsection{The non-coprime case}

Throughout this paper we have assumed for simplicity that the integers $s$ and $t$ are coprime. In fact, this assumption is unnecessary for many of our results. \cref{olssonnath} was generalised to the non-coprime case by Nath \cite[Theorem 1.1]{nath} and (with a different proof) by Gramain and Nath \cite[Theorem 2.1]{gn}; the idea in the latter proof is to consider the $g$-quotient of a partition, where $g$ is the greatest common divisor of $s$ and $t$. Applying \cref{olssonnath} to the components of this quotient and using results on quotients such as \cref{quotprops}, one obtains the general result. This technique can be applied to many of our results, too, and \cref{ineq,betacondition,stts,speccore,uniquest} all hold without modification in the case where $s$ and $t$ are not coprime, while \cref{basicprops} requires minor modification. The results concerning the level $t$ action of $\weyl s$ do not generalise so readily: one must consider the action of a group consisting of a direct product of $g$ copies of $\weyl{s/g}$. Then the results we have proved can be made to work, but the level of complication soon outweighs the reward; we leave the interested reader to work out the details. The results in \cref{kappasec} seem to have no analogue in the non-coprime case, where there are infinitely many $(s,t)$-cores.

\subsection{\gcc stus}

A natural extension of the results in this paper would be to try to extend from two integers $s,t$ to three (or more): is there is a suitable definition of an \gcc stu? The author has not been able to find the appropriate generalisation of our initial definition of an \gc st. However, \cref{uniquest} suggests a possibility: assuming $s,t,u>1$, we could define an \gcc stu to be a partition which is uniquely determined by its size and its $s$-, $t$- and $u$-cores. We hope to be able to say something about such partitions in the future.

\end{document}